\documentclass{amsart}
\usepackage{%
amssymb,%
amsrefs,
fullpage}

\input xypic
\xyoption{all}
\newdir{ >}{{}*!/-9pt/\dir{>}}

\newcommand{\Ab}        {\operatorname{Ab}}
\newcommand{\Aut}       {\operatorname{Aut}}
\newcommand{\End}       {\operatorname{End}}
\newcommand{\Ext}       {\operatorname{Ext}}
\newcommand{\Ho}        {\operatorname{Ho}}
\newcommand{\Hom}       {\operatorname{Hom}}
\newcommand{\Map}       {\operatorname{Map}}

\newcommand{\point}     {\operatorname{point}}

\newcommand{\CA}        {{\mathcal{A}}}
\newcommand{\CB}        {{\mathcal{B}}}
\newcommand{\CC}        {{\mathcal{C}}}
\newcommand{\CF}        {{\mathcal{F}}}
\newcommand{\CH}        {{\mathcal{H}}}
\newcommand{\CM}        {{\mathcal{M}}}
\newcommand{\CS}        {{\mathcal{S}}}
\newcommand{\CV}        {{\mathcal{V}}}

\newcommand{\N}         {{\mathbb{N}}}
\newcommand{\Z}         {{\mathbb{Z}}}
\newcommand{\Zm}[1]     {\mathbb{Z}/#1}              
\newcommand{\Fp}        {{\mathbb{F}_p}}          
\newcommand{\Q}         {{\mathbb{Q}}}
\newcommand{\R}         {{\mathbb{R}}}
\newcommand{\C}         {{\mathbb{C}}}

\newcommand{\Sg}        {\Sigma}
\newcommand{\Om}        {\Omega}

\newcommand{\al}        {\alpha}
\newcommand{\dl}        {\delta}
\newcommand{\ep}        {\epsilon}
\newcommand{\sg}        {\sigma}

\newcommand{\Sgi}       {\Sigma^\infty}
\newcommand{\Smash}     {\wedge}
\newcommand{\Wedge}     {\vee}
\newcommand{\bigWedge}  {\bigvee}
\newcommand{\hA}	{\widehat{A}}
\newcommand{\op}        {\oplus}
\newcommand{\ot}        {\otimes}
\newcommand{\smll}      {\text{small}}
\newcommand{\sse}       {\subseteq}
\newcommand{\st}        {\;|\;}
\newcommand{\tH}        {\widetilde{H}}
\newcommand{\tm}        {\times}
\newcommand{\xla}       {\xleftarrow}
\newcommand{\xra}       {\xrightarrow}
\newcommand{\psb}[1]    {[\![#1]\!]}
\newcommand{\ip}[1]     {\langle #1\rangle}

\newcommand{\colim}  {\operatornamewithlimits{\underset{\longrightarrow}{lim}}}
\newcommand{\invlim} {\operatornamewithlimits{\underset{\longleftarrow}{lim}}}

\renewcommand{\:}{\colon}

\theoremstyle{definition}

\begin{document}
\title{An introduction to the category of spectra}
\author{N.~P.~Strickland}
\bibliographystyle{abbrv}

\maketitle 

These notes give a brief introduction to the category of spectra as
defined in stable homotopy theory.  In particular,
Section~\ref{sec-examples} discusses an extensive list of examples of
spectra whose properties have been found to be interesting.  Although
many references are given, most of them are old.  A survey of more
recent research and exposition would be valuable, but is not attempted
here.

\section{Introduction}

Early in the history of homotopy theory, people noticed a number of
phenomena suggesting that it would be convenient to work in a context
where one could make sense of negative-dimensional spheres.  Let $X$
be a finite pointed simplicial complex; some of the relevant phenomena
are as follows.
\begin{itemize}
 \item For most $n$, the homotopy sets $\pi_nX$ are abelian groups.
  The proof involves consideration of $S^{n-2}$ and so breaks down for
  $n<2$; this would be corrected if we had negative spheres.
 \item Calculation of homology groups is made much easier by the
  existence of the suspension isomorphism $\tH_{n+k}\Sg^k X=\tH_nX$.
  This does not generally work for homotopy groups.  However, a
  theorem of Freudenthal says that if $X$ is a finite complex, we at
  least have a suspension isomorphism
  $\pi_{n+k}\Sg^k X=\pi_{n+k+1}\Sg^{k+1}X$ for large $k$.  If could
  work in a context where $S^{-k}$ makes sense, we could smash
  everything with $S^{-k}$ to get a suspension isomorphism in homotopy
  parallel to the one in homology.
 \item We can embed $X$ in $S^{k+1}$ for large $k$, and let $Y$ be the
  complement.  Alexander duality says that $\tH_nY=\tH^{k-n}X$,
  showing that $X$ can be ``turned upside-down'', in a suitable
  sense.  The shift by $k$ is unpleasant, because the choice of $k$ is
  not canonical, and the minimum possible $k$ depends on $X$.
  Moreover, it is unsatisfactory that the homotopy type of $Y$ is not
  determined by that of $X$ (even after taking account of $k$).  In a
  context where negative spheres exist, one can define
  $DX=S^{-k}\Smash Y$; one finds that $\tH_nDX=\tH^{-n}X$ and that
  $DX$ is a well-defined functor of $X$, in a suitable sense.
 \item The Bott periodicity theorem says that the homotopy groups of
  the infinite orthogonal group $O(\infty)$ satisfy
  $\pi_{k+8}O(\infty)=\pi_kO(\infty)$ for all $k\geq 0$.  It would be
  pleasant and natural to extend this pattern to negative values of
  $k$, which would again require negative spheres.
\end{itemize}

Considerations such as these led to the construction of the
Spanier-Whitehead category $\CF$ of finite spectra, which we briefly
survey in Section~\ref{sec-finite-spectra}.  Although fairly
straightforward, and very beautiful and interesting, this category has
two defects.
\begin{itemize}
 \item Many of the most important examples in homotopy theory are
  infinite complexes: Eilenberg-MacLane spaces, classifying spaces of
  finite groups, infinite-dimensional grassmannians and so on.  The
  category $\CF$ is strongly tied to finite complexes, so a wider
  framework is needed to capture these examples.
 \item Ordinary homotopy theory is made both easier and more
  interesting by its connections with geometry.  However, $\CF$ is
  essentially a homotopical category, with no geometric structure
  behind it.  This also prevents a good theory of spectra with a group
  action, or of bundles of spectra over a space, or of diagrams of
  spectra. 
\end{itemize}
The first problem was addressed by a number of people, but the
definitive answer was provided by Boardman.  He constructed a
category $\CB$ with excellent formal properties parallel to those of
$\CF$, whose subcategory of finite objects (suitably defined) is
equivalent to $\CF$.  A popular exposition of this category is in
Adams' book~\cite{ad:shg}.  Margolis~\cite{ma:ssa} gave a
list of the main formal properties of $\CB$ and its relationship with
$\CF$.  He conjectured (with good evidence) that they characterise
$\CB$ up to equivalence.  See~\cite{chst:pmh} for some new evidence
for this conjecture, and~\cite{hopast:ash} for an investigation of
some related systems of axioms.

The second problem took much longer to resolve.  There have been
a number of constructions of topological categories whose associated
homotopy category (suitably defined) is equivalent to $\CB$, with
steadily improving formal
properties~\cites{lemast:esh,el:ggs,ekmm:rma,hoshsm:ss}.  There is also
a theorem of Lewis~\cite{le:ccs} which shows that it is impossible to
have all the good properties that one might naively hope for.  We will
sketch one construction in Section~\ref{sec-boardman}.

\section{The finite stable category}
\label{sec-finite-spectra}

\subsection{Basics}

We first recall some basic definitions.  In this section all spaces
are assumed to be finite CW complexes with basepoints.  (We could
equally well use simplicial complexes instead, at the price of having
to subdivide and simplicially approximate from time to time.)  We
write $0$ for all basepoints, and we write $[A,B]$ for the set of
based homotopy classes of maps from $A$ to $B$.  We define $A\Wedge B$
to be the quotient of the disjoint union of $A$ and $B$ in which the
two basepoints are identified together.  We also define $A\Smash B$ to
be the quotient of $A\tm B$ in which the subspace $A\tm 0\cup 0\tm B$
is identified with the single point $(0,0)$.  This is called the smash
product of $A$ and $B$; note that it is commutative and associative up
to isomorphism and that $S^0\Smash A=A$, where $S^0=\{0,1\}$.
Moreover, $A\Smash(B\Wedge C)=(A\Smash B)\Wedge(A\Smash C)$ and
$[A\Wedge B,C]=[A,C]\tm[B,C]$ (so $\Wedge$ is the coproduct in the
homotopy category of pointed spaces).  We let $S^1$ denote the
quotient of $[0,1]$ in which $0$ is identified with $1$, and we write
$\Sg A=S^1\Smash A$.  One can check that $\Sg$ sends the ball $B^n$ to
$B^{n+1}$ and the sphere $S^n$ to $S^{n+1}$, and the reduced homology
of $\Sg A$ is just $\tH_n(\Sg A)=\tH_{n-1}(A)$.  Thus, we think of
$\Sg$ as shifting all dimensions by one.

The quotient space of $S^1$ in which $1/2$ is identified with $0$ is
evidently homeomorphic to $S^1\Wedge S^1$, so we get a map
$\dl\:S^1\xra{}S^1\Wedge S^1$ and thus a map 
$\dl\:\Sg A\xra{}\Sg A\Wedge\Sg A$.  It is well-known that the induced
map $\dl^*\:[\Sg A,B]\tm[\Sg A,B]\xra{}[\Sg A,B]$ makes $[\Sg A,B]$
into a group.  There are apparently $n$ different group structures on
$[\Sg^n A,B]$, but it is also well-known that they are all the same,
and they are commutative when $n>1$.  We have an evident sequence of
maps 
\[ [A,B]\xra{\Sg}[\Sg A,\Sg B]\xra{\Sg}[\Sg^2A,\Sg^2B]\xra{}\ldots. 
\]
Apart from the first two terms, it is a sequence of Abelian groups and
homomorphisms.  By a fundamental theorem of Freudenthal, after a
finite number of terms, it becomes a sequence of isomorphisms.  We
define $[\Sgi A,\Sgi B]$ to be the group $[\Sg^NA,\Sg^NB]$ for large
$N$, or if you prefer the colimit $\colim_N[\Sg^NA,\Sg^NB]$.  After
doing a little point-set topology, one concludes that this is the same
as the set $[A,QB]$, where $QB=\colim_N\Om^N\Sg^NB$ and $\Om^NC$ means
the space of based continuous maps $S^N\xra{}C$, with a suitable
topology.  

\subsection{Finite spectra}

One can define a category with one object called $\Sgi A$ for each
finite CW complex $A$, and morphisms $[\Sgi A,\Sgi B]$.  It is easy to
see that $\Sg$ induces a full and faithful endofunctor of this
category.  We prefer to arrange things so that $\Sg$ is actually an
equivalence of categories.  Accordingly, we define a category $\CF$
whose objects are expressions of the form $\Sg^{\infty+n}A$ where $A$
is a finite CW complex and $n$ is an integer.  (If you prefer, you can
take the objects to be pairs $(n,A)$.)  We refer to these objects as
finite spectra.  The maps are
\[ [\Sg^{\infty+n}A,\Sg^{\infty+m}B] =
    \colim_N[\Sg^{N+n}A,\Sg^{N+m}B].
\]
Freudenthal's theorem again assures us that the limit is attained at a
finite stage.  The functor $\Sg$ induces a self-equivalence of the
category $\CF$.  There are evident extensions of the functors $\Wedge$
and $\Smash$ to $\CF$ such that $\Sgi A\Wedge\Sgi B=\Sgi(A\Wedge B)$
and $\Sgi A\Smash\Sgi B=\Sgi(A\Smash B)$ (although care is needed with
signs when defining the smash product of morphisms).  The category
$\CF$ is additive, with biproduct given by the functor $\Wedge$.  The
morphism sets $[X,Y]$ in $\CF$ are finitely generated Abelian groups.
One can define homology of finite spectra by
$H_n\Sg^{\infty+m}A=\tH_{n-m}A$, and then the map
\[ H_*\: \Q\ot[X,Y] \xra{} \prod_n\Hom(H_n(X;\Q),H_n(Y;\Q)) \]
is an isomorphism.  The groups $[X,Y]$ themselves are known to be
recursively computable, but the guaranteed algorithms are of totally
infeasible complexity.  Nonetheless, there are methods of computation
which require more intelligence than the algorithms but have a
reasonable chance of success.  

\subsection{Stable homotopy groups of spheres}

Even the groups $\pi_n^S=[\Sg^{\infty+n}S^0,\Sgi S^0]$ are hard, and
are only known for $n\le 60$ or so (they are zero when $n<0$).  The
first few groups are as follows:
\begin{align*}
 \pi_0^S &= \Z\{\iota\}                                 &
 \pi_1^S &= \Zm{2}\{\eta\}                              &
 \pi_2^S &= \Zm{2}\{\eta^2\}                            &
 \pi_3^S &= \Zm{24}\{\nu\}                              \\
 \pi_4^S &= 0                                           &
 \pi_5^S &= 0                                           &
 \pi_6^S &= \Zm{2}\{\nu^2\}                             &
 \pi_7^S &= \Zm{240}\{\sg\}                             
\end{align*}
Here $\iota$ is the identity map, and $\eta$ comes from the map 
\[ \eta\:S^3=\{(z,w)\in\C^2\st |z|^2+|w|^2=1\}\xra{}
             \C\cup\{\infty\}=S^2
\]
defined by $\eta(z,w)=z/w$.  Similarly, $\nu$ comes from division of
quaternions, and $\sg$ from division of octonions (but most later
groups cannot be described in a similarly explicit way).  The
expression $\eta^2$ really means $\eta\circ(\Sg\eta)$, and $\nu^2$
means $\nu\circ(\Sg^3\nu)$.

Many general results are also known.  For example, for any prime $p$,
the $p$-torsion part of $\pi_n^S$ is known for $n<2p^3-2p$ and is zero
for $n<2p-3$ (provided $n\neq 0$).  Both the rank and the exponent are
finite but unbounded as $n$ tends to infinity.  The group $\pi_*^S$ is
a graded ring, and is commutative in the graded sense.  An important
theorem of Nishida says that all elements of degree greater than zero
are nilpotent.

\subsection{Triangulation}

The category $\CF$ is not Abelian.  Instead, it has a triangulated
structure.  This means that there is a distinguished class of diagrams
of the shape $X\xra{f}Y\xra{g}Z\xra{h}\Sg X$ (called exact triangles)
with certain properties to be listed below.  In our case the exact
triangles can be described as follows.  Let $A$ be a subcomplex of a
finite CW complex $B$, and let $C$ be obtained from $B$ by attaching a
cone $I\Smash A$ along the subspace $\{1\}\tm A=A$.  There is an
evident copy of $B$ in $C$, and if we collapse it to a point we get a
copy of $\Sg A$.  We thus have a diagram of spaces
$A\xra{}B\xra{}C\xra{}\Sg A$.  We say that a diagram
$X\xra{}Y\xra{}Z\xra{}\Sg X$ of finite spectra is an exact triangle if
it is isomorphic to a diagram of the form
$\Sg^{\infty+n}A\xra{}\Sg^{\infty+n}B\xra{}
 \Sg^{\infty+n}C\xra{}\Sg^{\infty+n+1}A$ for some $n\in\Z$ and some
$A$, $B$ and $C$ as above.  Incidentally, one can show that $C$ is
homotopy equivalent to the space $B/A$ obtained from $B$ by
identifying $A$ with the basepoint.

The axioms for a triangulated category are as follows.  In our case,
they all follow from the theory of Puppe sequences in unstable
homotopy theory.
\begin{itemize}
 \item[(a)] Any diagram isomorphic to an exact triangle is an exact
  triangle. 
 \item[(b)] Any diagram of the form $0\xra{}X\xra{1}X\xra{}\Sg 0=0$ is
  an exact triangle.
 \item[(c)] Any diagram $X\xra{f}Y\xra{g}Z\xra{h}\Sg X$ is an exact
  triangle if and only if the diagram
  $Y\xra{g}Z\xra{h}\Sg X\xra{-\Sg f}\Sg X$ is an exact triangle.
 \item[(d)] For any map $f\:X\xra{}Y$ there exists a spectrum $Z$ and
  maps $g,h$ such that $X\xra{f}Y\xra{g}Y\xra{h}\Sg X$ is an exact
  triangle.
 \item[(e)] Suppose we have a diagram as shown below (with $h$
  missing), in which the rows are exact triangles and the
  rectangles commute.  Then there exists a (nonunique) map $h$ making
  the whole diagram commutative.
  \[ \xymatrix{
      U \ar[r] \ar[d]_f &
      V \ar[r] \ar[d]_g &
      W \ar[r] \ar@{.>}[d]^h &
      \Sg U \ar[d]^{\Sg f} \\
      X \ar[r] &
      Y \ar[r] &
      Z \ar[r] &
      \Sg X
  } \]
 \item[(f)] Suppose we have maps $X\xra{v}Y\xra{u}Z$, and exact
  triangles $(X,Y,U)$, $(X,Z,V)$ and $(Y,Z,W)$ as shown in the
  diagram.  (A circled arrow
  $U\longrightarrow\hspace{-1.3em}\circ\hspace{0.8em} X$ means a map
  $U\xra{}\Sigma X$.)  Then there exist maps $r$ and $s$ as shown,
  making $(U,V,W)$ into an exact triangle, such that the following
  commutativities hold:
  \[ au = rd \qquad \qquad es = (\Sigma v) b \qquad\qquad
     sa = f  \qquad \qquad br = c \]
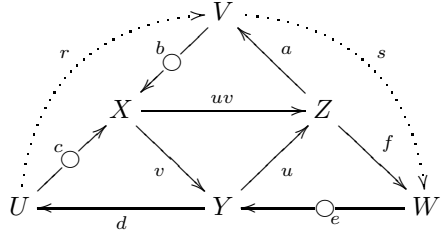
  \[ \xymatrix{
      & & V  \ar@{.>}@/^2pc/[ddrr]^s \ar@{->}[dl]_b|\bigcirc \\
      & X \rrto^{uv} \drto_v & &
      Z \ulto_a \drto^f \\
      U  \ar@{.>}@/^2pc/[uurr]^r \urto^c|\bigcirc & &
      Y \llto^d \urto_u  & &
      W \llto^(0.45)e|\bigcirc
     }
  \]
\end{itemize}
The last axiom is called the octahedral axiom (the diagram can be
turned into an octahedron by lifting the outer vertices and drawing an
extra line from $W$ to $U$).  In our case it basically just says that
when we have inclusions $A\sse B\sse C$ of CW complexes we have
$(C/A)/(B/A)=C/B$.  

One of the most important consequences of the axioms is that whenever
$X\xra{}Y\xra{}Z\xra{}\Sg X$ is an exact triangle and $W$ is a finite
spectrum, we have long exact sequences
\[ \ldots \xra{} [W,\Sg^{-1}Z]\xra{} [W,X] \xra{}
       [W,Y] \xra{} [W,Z] \xra{} [W,\Sg X] \xra{}\ldots
\]
and
\[ \ldots \xla{} [\Sg^{-1}Z,W] \xla{} [X,W] \xla{}
       [Y,W] \xla{} [Z,W] \xla{} [\Sg X,W] \xla{} \ldots.
\]

\subsection{Thom spectra}

Let $X$ be a finite CW complex, and let $V$ be a vector bundle over
$X$.  The Thom space $X^V$ can be defined as the one-point
compactification of the total space of $V$.  This has many interesting
properties, not least of which is the fact that when $V$ is an
oriented bundle of dimension $n$, the reduced cohomology $\tH^*(X^V)$
is a free module over $H^*(X)$ on one generator in dimension $n$.
This construction can be generalised to virtual bundles, in other
words formal expressions of the form $V-W$, except that we now have a
Thom spectrum $X^{V-W}$ rather than a Thom space.  The construction is
to choose a surjective map from a trivial bundle $\R^n\tm X$ onto
$W$, with kernel $U$ say, and define 
$X^{V-W}=\Sg^{\infty-n}X^{V\op U}$.

\subsection{Duality}
\label{subsec-duality}

For any finite spectrum $X$, there is an essentially unique spectrum
$DX$ (called the Spanier-Whitehead dual of $X$) equipped with a
natural isomorphism $[W\Smash X,Y]=[W,DX\Smash Y]$.  This can be
constructed in a number of different ways.  One way is to start with a
simplicial complex $A$ and embed it simplicially as a proper
subcomplex of $S^{N+1}$ for some $N>0$.  One can show that the
complement of $A$ has a deformation retract $B$ which is a finite
simplicial complex, and $D\Sg^{\infty+n}A=\Sg^{\infty-n-N}B$.  Note
that Alexander duality implies that $H_mX=H^{-m}DX$.  

An important example arises when $X=\Sgi M_+$ for some smooth manifold
$M$, with tangent bundle $\tau$ say.  It is not hard to show
geometrically that $D(\Sgi M_+)$ is $M^{-\tau}$, the Thom spectrum of
the virtual bundle $-\tau$ over $M$; this phenomenon is called
\emph{Atiyah duality}.

We also write $F(X,Y)=DX\Smash Y$.  This is a functor in both
variables, it preserves cofibrations up to sign, and the defining
property of $DX$ can be rewritten as $[W,F(X,Y)]=[W\Smash X,Y]$.

\subsection{Splittings}
\label{subsec-splittings}

It often happens that we have a finite complex $X$ that cannot be
split into simpler pieces, but that the finite spectrum $\Sgi X$ does
have a splitting.  Group actions are one fruitful source of
splittings.  If a finite group $G$ acts on $X$, then the map
$G\xra{}\Aut(\Sgi X)$ extends to a ring map $\Z[G]\xra{}\End(\Sgi X)$.
If $\tH_*X$ is a $p$-torsion group, then this will factor through
$(\Z/p^n)[G]$ for large $n$.  Any idempotent element in $(\Z/p)[G]$
can be lifted uniquely to an idempotent in $(\Z/p^n)[G]$, which will
give an idempotent in $\End(\Sgi X)$ and thus a splitting of $X$.  The
methods of modular representation theory give good information about
idempotents in group rings, and thus a supply of interesting
splittings.  The Steinberg idempotent in $(\Z/p)[GL_n(\Z/p)]$ gives
particularly important examples~\cite{mi:osm,mipr:ssd}, as do various
idempotents in $(\Z/p)[\Sg_n]$ (see~\cite{ra:nps}*{Appendix C}).

Another common situation is to have a finite complex $X$ and a
filtration $F_0X\sse F_1X\sse \ldots\sse F_nX=X$ that splits stably,
giving an equivalence $\Sgi X\simeq\bigWedge_n\Sgi F_nX/F_{n-1}X$ of
finite spectra.  For example, one can take $X=U(m)$, and let $F_nX$ be
the space of matrices $A\in U(m)$ for which the rank of $A-I$ is at
most $n$.  A theorem of Miller~\cite{mi:sss} says that the filtration
splits stably, and that the quotient $F_nU(m)/F_{n-1}U(m)$ is the Thom
space of a certain bundle over the Grassmannian of $n$-planes in
$\C^m$.  Later we will explain how to interpret $\Sgi X$ when $X$ is
an infinite complex; there are many examples in which $X$ has a stably
split filtration in which the quotients are finite spectra.  This
holds for $X=BU(n)$ or $X=\Om U(n)$ or $X=\Om^nS^{n+m}$, for example.
The splitting of $\Om^nS^{n+m}$ is due to Snaith~\cite{sn:sdo}; the
Snaith summands in $\Om^2S^3$ are called \emph{Brown-Gitler spectra},
and they have interesting homological properties with many
applications~\cite{brgi:scc,sc:ums,huku:csu}.

\section{Cobordism and Morava $K$-theory}
\label{sec-cobordism}

We next outline the theory of complex cobordism~\cite{ra:ccs}, and the
results of Hopkins, Devinatz and Smith showing how complex cobordism
reveals an important part of the structure of $\CF$.

Given a space $X$ and an integer $n\geq 0$ we define a \emph{geometric
  $n$-chain} in $X$ to be a compact smooth $n$-manifold $M$ (possibly with
boundary) equipped with a continuous map $f\:M\xra{}X$.  We regard
$(M_0,f_0)$ and $(M_1,f_1)$ as equivalent if there is a
differomorphism $g\:M_0\xra{}M_1$ with $f_1g=f_0$.  We write $GC_nX$
for the set of equivalence classes, which is a commutative semigroup
under disjoint union.  We define a differential
$\partial\:GC_nX\xra{}GC_{n-1}X$ by $\partial[M,f]=[\partial
M,f|_{\partial M}]$.  One can make sense of the homology
$MO_*X=H(GC_*X,\partial)$, and (because $\partial(M\tm I)=M\amalg M$)
one finds that $MO_*X$ is a vector space over $\Z/2$.  In the case
where $X$ is a point, one can use cartesian products to make
$MO_*(\point)$ into a graded ring, which is completely described by a
remarkable theorem of Thom: it is a polynomial algebra over $\Z/2$
with one generator $x_n$ in degree $n$ for each integer $n>0$ not of
the form $2^k-1$.  New perspectives on this answer and the underlying
algebra were provided by Quillen~\cite{qu:fgl,qu:epc} and
Mitchell~\cite{mi:psm}.  One can also show that
$MO_*X=MO_*(\point)\ot_{\Z/2}H_*(X;\Z/2)$, so this construction does
not yield new invariants of spaces.  This isomorphism gives an obvious
way to define $MO_*X$ when $X$ is a finite spectrum.

The story changes however, if we work with oriented manifolds.  This
gives groups $MSO_*X$ with a richer structure; in particular, they are
not annihilated by $2$.  There are various hints that complex
manifolds would give still more interesting invariants, but there are
technical problems, not least the lack of a good theory of complex
manifolds with boundary.  It turns out to be appropriate to generalize
and consider manifolds with a specified complex structure on the
stable normal bundle, known as ``stably complex manifolds''.  The
precise definitions are delicate; details are explained
in~\cite{qu:epc}, and Buchstaber and Ray~\cite{bura:tst} have provided
naturally occuring examples where the details are important.  In any
case, one ends up with a ring $MU_*=MU_*(\point)$, and groups $MU_*X$
for all spaces $X$ that are modules over it.  It is not the case that
$MU_*X=MU_*\ot_\Z H_*X$, but there is still a suspension isomorphism,
which allows one to define $MU_*X$ when $X$ is a finite spectrum.  One
finds that this is a generalized homology theory (known as
\emph{complex cobordism}), so it converts cofibre seqences of spectra
to long exact sequences of modules.  The nilpotence theory of
Devinatz, Hopkins and Smith~\cite{dehosm:nsh,hosm:nshii,ra:nps} (which
will be outlined below) shows that $MU_*X$ is an extremely powerful
invariant of $X$.
 
The ring $MU_*$ turns out to be a polynomial algebra over $\Z$ with
one generator in each positive even degree.  There is no canonical
system of generators, but nonetheless, Quillen showed that $MU_*$ is
canonically isomorphic to an algebraically defined object: Lazard's
classifying ring for formal group laws~\cites{ad:shg,ra:ccs}.  This
was the start of an extensive relationship between stable homotopy and
formal group theory.  The algebra provides many natural examples of
graded rings $A_*$ equipped with a formal group law and thus a map
$MU_*\xra{}A_*$.  It is natural to ask whether there is a generalised
homology theory $A_*(X)$ whose value on a point is the ring $A_*$.
Satisfactory answers for a broad class of rings $A_*$ are given
in~\cite{st:pmm}, which surveys and consolidates a great deal of much
older literature and extends newer ideas from~\cite{ekmm:rma}.  In
particular, we can consider the rings $K(p,n)_*=\Fp[v_n,v_n^{-1}]$,
where $p$ is prime and $n>0$ and $v_n$ has degree $2(p^n-1)$.  This is
made into an algebra over $MU_*$ using a well-known formal group law,
and by old or new methods, one can construct an associated generalised
homology theory $K(p,n)_*X$, known as \emph{Morava $K$-theory}.  It is
convenient to extend the definition by putting $K(p,0)_*X=H_*(X;\Q)$.

A key theorem of Hopkins, Devinatz and Smith says that if
$f\:\Sg^dX\xra{}X$ is a self-map of a finite spectrum $X$, and
$K(p,n)_*(f)=0$ for all primes $p$ and all $n\geq 0$, then the
iterated composite $f^m\:\Sg^{md}X\xra{}X$ is zero for large $m$, or
in other words $f$ is composition-nilpotent.  They also show that
Morava $K$-theory detects nilpotence in a number of other senses, and
give formulations involving the single theory $MU$ instead of the
collection of theories $K(p,n)$.  

If $R$ is a commutative ring, it is well-known that the ideal of
nilpotent elements is the intersection of all the prime ideals, so the
Zariski spectrum is unchanged if we take the quotient by this ideal.
One can deduce that the classification of certain types of
subcategories of the abelian category of $R$-modules is again
insensitive to the ideal of nilpotents.  Our category $\CF$ of finite
spectra is triangulated rather than abelian, but nonetheless Hopkins
and Smith developed an analogous theory and deduced a classification
of the thick subcategories of $\CF$, with many important
consequences. 

\section{Boardman's category $\CB$}
\label{sec-boardman}

It is clearly desirable to have a category $\CC$ analogous to $\CF$
but without finiteness conditions.  There are various obvious
candidates: one could take the Ind-completion of $\CF$, or just follow
the definition of $\CF$ but allow infinite CW complexes instead of
finite ones.  Unfortunately, these categories turn out to have
unsatisfactory technical properties.  The requirements were first
assembled in axiomatic form by Margolis~\cite{ma:ssa}; in outline,
they are as follows:
\begin{itemize}
 \item $\CC$ should be a triangulated category
 \item Every family $\{X_\al\}$ of objects in $\CC$ should have a
  coproduct, written $\bigWedge_\al X_\al$
 \item For any $X,Y\in\CC$ there should be functorially associated
  objects $X\Smash Y$ and $F(X,Y)$ making $\CC$ a closed symmetric
  monoidal category.
 \item If we let $\text{small}(\CC)$ be the subcategory of objects $W$
  for which the natural map
  $\bigoplus_\al[W,X_\al]\xra{}[W,\bigWedge_\al X_\al]$ is always an
  isomorphism, then $\smll(\CC)$ should be equivalent to $\CF$.
\end{itemize}
Historically, the work of Margolis came after Boardman's construction
of a category $\CB$ satisfying the axioms, and Adams's
explanation~\cite{ad:shg} of a slightly different way to approach the
construction.  Margolis conjectured that if $\CC$ satisfies the axioms
then $\CC$ is equivalent to $\CB$.  Schwede and
Shipley~\cite{scsh:uts} have proved that this is true, provided that
$\CC$ is the homotopy category of a closed model category in the sense
of Quillen~\cite{qu:ha,dwsp:htm} satisfying suitable axioms.  There is
also good evidence for the conjecture without this additional
assumption.  The objects of $\CB$ are generally called \emph{spectra},
although in some contexts one introduces different words to
distinguish between objects in different underlying geometric
categories. 

Probably the best approach to constructing $\CB$ is via the theory of
\emph{orthogonal spectra}, as we now describe~\cite{mamash:mcd}.  Let
$\CV$ denote the category of finite-dimensional vector spaces over
$\R$ equipped with an inner product.  The morphisms are linear
isomorphisms that preserve inner products.  For any $V\in\CV$ (with
$\dim(V)=n$ say) we write $S^V$ for the one-point compactification of
$V$; this is homeomorphic to $S^n$.  An orthogonal spectrum $X$
consists of a functor $\CV\to\{\text{based spaces}\}$ together with
maps $S^U\Smash X(V)\to X(U\oplus V)$ satisfying various continuity
and compatibility conditions that we will not spell out.  We write
$\CS$ for the category of orthogonal spectra.  For example, a based
space $A$ gives an orthogonal spectrum $\Sgi A$ with
$(\Sgi A)(U)=S^U\Smash A$.  A vector space $T\in\CV$ gives positive
and negative sphere spectra $S^{\pm T}$: the value of the sphere
spectrum $S^T$ at $U\in\CV$ is the sphere space
$S^U\Smash S^T=S^{U\oplus T}$, and
\[ S^{-T}(U) =
    \{(\al,u)\st \al\in\CV(T,U),\; u\in U,\;\ip{u,\al(T)}=0\}
    \cup \{\infty\}.
\]
   
For orthogonal spectra $X$ and $Y$, the morphism set $\CS(X,Y)$ has a
natural topology, and we could define an associated homotopy category
by the rule $[X,Y]=\pi_0\CS(X,Y)$.  Unfortunately, the resulting
category is not the one that we want.  Instead, we define the homotopy
groups of $X$ by the rule $\pi_k(X)=\colim_N\pi_{k+N}(X(\R^N))$.  We
then say that a map $f\:X\to Y$ is a \emph{weak equivalence} if
$\pi_*(f)\:\pi_*(X)\to\pi_*(Y)$ is an isomorphism.  We now construct a
new category $\Ho(\CS)$ by starting with $\CS$ and adjoining formal
inverses for all weak equivalences.  It can be shown that this is
equivalent to $\CB$ (or can be taken as the definition of $\CB$).

This process of adjoining formal inverses can be subtle.  To manage
the subtleties, we need the theory of model categories in the sense of
Quillen~\cites{dwsp:htm,ho:mc}.  In particular, this will show that
$\Ho(\CS)(X,Y)=\pi_0\CS(X,Y)$ for certain classes of spectra $X$ and
$Y$; this is enough to get started with computations and prove that
Margolis's axioms are satisfied.

Given orthogonal spectra $X$, $Y$ and $Z$, a \emph{pairing} from $X$
and $Y$ to $Z$ consists of maps $\al_{U,V}\:X(U)\Smash Y(V)\to
Z(U\oplus V)$ satisfying some obvious compatibility conditions.  One
can show that there is an orthogonal spectrum $X\Smash Y$ such that
pairings from $X$ and $Y$ to $Z$ biject with morphisms from
$X\Smash Y$ to $Z$.  Basic examples are that $S^0\Smash X=X$ and
$\Sgi A\Smash\Sgi B=\Sgi(A\Smash B)$ and $S^U\Smash S^V=S^{U\oplus V}$
and $S^{-U}\Smash S^{-V}=S^{-U-V}$.  There is a natural map
$S^U\Smash S^{-U}\to S^0$ that is a weak equivalence but not an
isomorphism. 

This construction gives a symmetric monoidal structure on
$\CS$.  This in turn gives rise to a symmetric monoidal structure on
$\Ho(\CS)$; however, there are some hidden subtleties in this step,
which again are best handled by the general theory of model
categories.  One consequence is that the topology of the classifying
space of the symmetric group $\Sg_k$ is mixed in to the structure of
the $k$-fold smash product $X^{(k)}=X\Smash\dotsb\Smash X$ and the
quotient $X^{(k)}/\Sg_k$.  

A \emph{ring spectrum} is an object $R\in\Ho(\CS)$ equipped with a
unit map $\eta\:S^0\to R$ and a multiplication map 
$\mu\:R\Smash R\to R$ such that the following diagrams in $\Ho(\CS)$
commute:
\[ \xymatrix{
   R\Smash R\Smash R \ar[r]^-{\mu\Smash 1} \ar[d]_-{1\Smash\mu} &
   R\Smash R \ar[d]^-\mu & 
   R \ar[r]^-{\eta\Smash 1} \ar[dr]_-1 &
   R\Smash R \ar[d]^-\mu &
   R \ar[l]_-{1\Smash\eta} \ar[dl]^-1
   \\
   R\Smash R \ar[r]_-\mu &
   R & & R
} \]
Because we now have a good underlying geometric category $\CS$, we can
formulate a more precise notion: a \emph{strict ring spectrum} is an
object $R\in\CS$ equipped with morphisms
$S^0\xra{\eta}R\xla{\mu}R\Smash R$ such that the above diagrams
commute in $\CS$ (not just in $\Ho(\CS)$).  

The symmetric monoidal structure on $\CS$ includes a natural map
$\tau_{XY}\:X\Smash Y\to Y\Smash X$.  We say that a strict ring
spectrum $R$ is \emph{strictly commutative} if
$\mu\circ\tau_{RR}=\mu$.  This is a surprisingly stringent condition,
with extensive computational consequences.  A key point is that we
have an iterated multiplication map $R^{(k)}/\Sg_k\to R$, and this
brings into play the structure of $R^*(B\Sg_k)$.

\section{Examples of spectra}
\label{sec-examples}

Some important functors that construct objects of $\CB$ are as
follows:
\begin{itemize}
 \item[(a)] For any based space $X$ there is a suspension spectrum
  $\Sgi X\in\CB$, whose homotopy groups are given by
  $\pi_n\Sgi X=\pi^S_nX=\colim_k\pi_{n+k}\Sg^kX$.  The relevant
  orthogonal spectrum is just $(\Sgi X)(V)=S^V\Smash X$.

  We will mention one important example of infinite complexes $X$ and
  $Y$ for which $[\Sgi X,\Sgi Y]$ is well-understood.  Let $G$ be a
  finite group, with classifying space $BG$.  Let $AG^+$ be the set of
  isomorphism classes of finite sets with a $G$-action.  We can define
  addition and multiplication on $AG^+$ by $[X]+[Y]=[X\amalg Y]$ and
  $[X][Y]=[X\tm Y]$.  There are no additive inverses, but we can
  formally adjoin them to get a ring called $AG$, the Burnside ring of
  $G$; this is not hard to work with explicitly.  There is a ring map
  $\ep\:AG\xra{}\Z$ defined by $\ep([X]-[Y])=|X|-|Y|$, with kernel $I$
  say.  We then have a completed ring $\hA G=\invlim_nAG/I^n$.  The
  Segal conjecture (which was proved by Carlsson~\cite{ca:esh}) gives
  an isomorphism $\hA G\simeq[\Sgi BG_+,\Sgi S^0]$.  One can deduce a
  description of $[\Sgi BG_+,\Sgi BH_+]$ in similar terms for any
  finite group $H$.
 \item[(b)] For any virtual vector bundle $V$ over any space $X$,
  there is a Thom spectrum $X^V\in\CB$.  In particular, if $V$ is the
  tautological virtual bundle (of virtual dimension zero) over the
  classifying space $BU$, then there is an associated Thom spectrum,
  normally denoted by $MU$.  This has the property that the groups
  $MU_*X$ (as in Section~\ref{sec-cobordism}) are given by
  $\pi_*(MU\Smash\Sgi X_+)$ (this is proved by a geometric argument,
  and is essentially the first step in the calculation of $MU_*$).
  One can construct $MU$ (and also $MO$ and $MSO$) as strictly
  commutative ring spectra.

  Bott periodicity gives an equivalence $BU\simeq\Om SU$.  A range of
  different proofs can be found in
  \citelist{
   \cite{mi:mt}*{Section 23}
   \cite{prse:lg}*{Sections 6.4 and 8.8}
   \cite{hu:fb}*{Chapter 11}
   \cite{be:npb}}.
  The filtration of $\Om SU$ by the subspaces $\Om SU(k)$ gives a
  filtration of $MU$ by subspectra $X(m)$, which are important in the
  proof of the Hopkins-Devinatz-Smith nilpotence theorem.  There are
  models of these homotopy types that are strict ring spectra, but
  they cannot be made commutative.
 \item[(c)] For any generalized cohomology theory $A^*$, there is an
  essentially unique spectrum $A\in\CB$ with $A^nX=[\Sgi X_+,\Sg^nA]$
  for all spaces $X$ and $n\in\Z$.  Similarly, for any generalized
  homology theory $B_*$, there is an essentially unique spectrum
  $B\in\CB$ with $B_nX=\pi_n(B\Smash\Sgi X_+)$ for all spaces $X$ and
  $n\in\Z$.  These facts are known as \emph{Brown representability};
  the word ``essentially'' hides some subtleties.  
 \item[(d)] In particular, for any abelian group $A$ there is an
  Eilenberg-MacLane spectrum $HA\in\CB$ such that
  $[\Sgi X_+,\Sg^nHA]=H^n(X;A)$ and
  $\pi_n(HA\Smash\Sgi X_+)=H_n(X;A)$.  (In fact, Brown
  representability is not needed for this: we have an explicit
  orthogonal spectram $HA$ where $HA(V)$ is the tensor product of $A$
  with the free abelian group generated by $V$, equipped with a
  suitable topology.)  If $A$ is a commutative ring,
  then $HA$ is a strictly commutative ring spectrum.  It is common to
  consider the case $A=\Z/2$.  Here it can be shown that 
  \[ \pi_*((H\Z/2)\Smash(H\Z/2)) = 
      \Z/2[\xi_1,\xi_2,\xi_3,\dotsc],
  \]
  with $|\xi_k|=2^k-1$.  This is known as the \emph{dual Steenrod
   algebra}, and denoted by $\CA_*$.  The dual group
  $\CA^k=\Hom(\CA_k,\Z/2)$ can be identified with
  $[H\Z/2,\Sg^k H\Z/2]$, so these groups again form a graded ring
  (under composition), called the \emph{Steenrod algebra}.  This ring
  is noncommutative, but its structure can be described quite
  explicitly.   It is important, because the mod $2$ cohomology of any
  space (or spectrum) has a natural structure as an $\CA^*$-module.
  There is a similar story for mod $p$ cohomology when $p$ is an odd
  prime, but the details are a little more complicated.
 \item[(e)] Another consequence of Brown representability is that
  there is a spectrum $I\in\CB$ such that
  $[X,I]\simeq\Hom(\pi_0X,\Q/\Z)$ for all $X\in\CB$.  This is called
  the \emph{Brown-Comenetz dual}~\cite{brco:pdg} of $S^0$; it is
  geometrically mysterious, and a fertile source of counterexamples.
 \item[(f)] If $M_*$ is a flat module over $MU_*$, then the functor
  $X\mapsto M_*\ot_{MU_*}MU_*X$ is a homology theory, so there is a
  representing spectrum $M$ with
  $\pi_*(M\Smash X)=M_*\ot_{MU_*}MU_*X$.  The Landweber exact functor
  theorem shows that flatness is not actually necessary: a weaker
  condition called \emph{Landweber exactness} will
  suffice~\cite{la:hpc}.  This condition is formulated in terms of
  formal group theory, and is often easy to check in practice.

  Often $M_*$ is a ring, and the $MU_*$-module structure arises from a
  ring map $MU_*\to M_*$, which corresponds (by Quillen's description
  of $MU_*$) to a formal group law over $M_*$.  

  Important examples include the Johnson-Wilson spectra $E(p,n)$, with
  $E(p,n)_*=\Z_{(p)}[v_1,\dotsc,v_n][v_n^{-1}]$ (where
  $|v_k|=2(p^k-1)$).  This has a canonical formal group law, which we
  will not describe here.  The ring $K(p,n)_*=\Z/p[v_n,v_n^{-1}]$ is
  naturally a quotient $E(p,n)_*/I_n$, where
  $I_n=(p,v_1,\dotsc,v_{n-1})$.  It it is not Landweber exact, but a
  corresponding spectrum $K(p,n)$ can be constructed by other means.
  It is also useful to consider the completed spectra
  $\widehat{E}(p,n)$ with 
  \[ \pi_*(\widehat{E}(p,n)) = 
      (E(p,n)_*)^\wedge_{I_n} =
        \Z_p\psb{v_1,\dotsc,v_{n-1}}[v_n^{\pm 1}].
  \]
  This is again Landweber exact.

  When $M_*$ is a ring, one might hope to find a model of this
  homotopy type that is actually a strict ring spectrum, preferably
  strictly commutative.  Unfortunately, this does not work very well.
  Often there will be uncountably many different ways to make $M$ into
  a strict ring spectrum, with no way to pick out a preferred choice.
  Moreover, there will often not be any choice that is strictly
  commutative.  However, by a theorem of Hopkins and Miller, there is
  an essentially unique strictly commutative model for
  $\widehat{E}(p,n)$.  The reason why this case is special involves
  quite deep aspects of the algebraic theory of formal groups.

  If $M_*$ is an algebra over $MU_*$ that does not satisfy the
  Landweber criterion, one can still try to produce a corresponding
  spectrum $M$ by other methods.  As mentioned previously, the
  paper~\cite{st:pmm} contains results in this direction.
 \item[(g)] For any small symmetric monoidal category category $\CA$,
  there is a $K$-theory spectrum $K(\CA)\in\CB$.  Computationally,
  this is very mysterious, apart from the fact that it is always
  connective (ie $\pi_nK(\CA)=0$ for $n<0$) and $\pi_0K(\CA)$ is the
  group completion of the monoid of connected components in $\CA$.
  Thomason has shown~\cite{th:smc} that for every connective spectrum
  $X$ there exists $\CA$ with $K(\CA)\simeq X$.
  \begin{itemize}
   \item If $\CA$ is the category of finite sets and isomorphisms,
    then $K(\CA)=\Sgi S^0$. 
   \item Let $G$ be a finite group, and let $\CA$ be the category of
    finite $G$-sets and isomorphisms.  Let $\CA_f$ (resp $\CA_t$) be
    the subcategory of free (resp. transitive) $G$-sets.  Then
    $K(\CA)=\Sgi B(\CA_t)_+$, which can also be described as the wedge
    over the conjugacy classes of subgroups $H\leq G$ of the spectra
    $\Sgi BW_GH_+$, or as the fixed point spectrum of the
    $G$-equivariant sphere spectrum in the sense of
    Lewis-May-Steinberger~\cite{lemast:esh}.  On the other hand,
    $K(\CA_f)=\Sgi BG_+$.
   \item Work of Kathryn Lesh~\cite{le:fsa} can be interpreted as
    exhibiting symmetric monoidal categories $\CM_n$ (of ``finite
    multisets with multiplicities at most $n$'') whose $K$-theory is
    the $n$'th symmetric power of $\Sgi S^0$.
   \item If $\CA$ is the symmetric monoidal category with object set
    $\N$ and only identity morphisms, then $K(\CA)=H\Z$.
   \item One can set up a category $\CA$, whose objects are smooth
    compact closed $1$-manifolds, and whose morphisms are cobordisms
    between them.  With the right choice of details, the $K$-theory
    spectrum $K(\CA)$ is then closely related to the classifying space
    of the stable mapping class group, and an important theorem of
    Madsen and Weiss~\cite{mawe:sms} can be interpreted as saying that
    $K(\CA)$ is the Thom spectrum of the negative of the tautological
    bundle over $\C P^\infty$, up to adjustment of $\pi_{-2}$.  
   \item If $R$ is a commutative ring and $\CA$ is the category of
    finitely generated projective $R$-modules, then $K(\CA)$ is the
    algebraic $K$-theory spectrum usually denoted by $K(R)$.  Even in
    the case $R=\Z$, this contains a great deal of arithmetic
    information.  By rather different methods one can construct
    spectra called $THH(R)$ and $TC(R)$ (topological Hochschild
    homology and topological cyclic homology) that approximate $K(R)$;
    there is an extensive literature on these approximations.  The
    definitions can be set up in such a way that the spectra $K(R)$,
    $THH(R)$ and $TC(R)$ are all strictly commutative ring spectra.
  \end{itemize}
 \item[(h)] The above construction can be modified slightly to take
  account of a topology on the morphism sets of $\CA$.  We can then
  feed in the category of finite-dimensional complex vector spaces and
  isomorphisms (or a skeleton thereof) to get a spectrum known as
  $kU$, the connective complex $K$-theory spectrum, with a homotopy
  element $u\in\pi_2kU$ such that $\pi_*kU=\Z[u]$.  This has the
  property that for finite complexes $X$, the group
  $kU^0X=[\Sgi X_+,kU]$ is the group completion of the monoid of
  isomorphism classes of complex vector bundles on $X$, or in other
  words the $K$-theory of $X$ as defined by Atiyah and
  Hirzebruch~\cite{at:kt} (inspired by Grothendieck's similar
  definition in the context of algebraic geometry).  By a colimit
  construction, one can build a periodized version called $KU$ with
  $\pi_*KU=\Z[u,u^{-1}]$.  There are direct constructions of $kU$ and
  $KU$ using Bott periodicity rather than symmetric monoidal
  categories.  There are also constructions with a more analytic
  flavour, based on spaces of Fredholm operators and so
  on~\cite{jo:srs}.  With the right construction, both $kU$ and $KU$
  are strictly commutative ring spectra.  It can be shown that $KU$ is
  naturally a Landweber exact $MU$-algebra, so
  $KU_*(X)=KU_*\ot_{MU_*}MU_*(X)$ for all $X$.

  The infinite complex projective space $\C P^\infty$ is well-known to
  be a commutative group up to homotopy.  Using this, one can make the
  spectrum $R=\Sgi(\C P^\infty)_+$ into a ring spectrum.  The standard
  identification $\C P^1=S^2$ gives rise to an element $v\in\pi_2(R)$,
  and we can use a homotopy colimit construction to invert this
  element, giving a new ring spectrum $R[v^{-1}]$.  It is a theorem of
  Snaith~\cite{sn:lsh} that $R[v^{-1}]$ is homotopy equivalent to $KU$.
 \item[(i)] Let $F$ be a functor from based spaces to based spaces.
  Under mild conditions, we can use the homeomorphism
  $S^1\Smash S^n\xra{}S^{n+1}$ to get a map 
  \[ S^1 \xra{} \Map(S^n,S^{n+1}) \xra{F} \Map(FS^n,FS^{n+1}), \] and
  thus an adjoint map $\Sg FS^n\xra{}FS^{n+1}$.  This gives a sequence
  of spectra $\Sg^{-n}\Sgi FS^n$, whose homotopy colimit (in a
  suitable sense) is denoted by $D_1F$.  This is called the
  \emph{linearization} or \emph{first Goodwillie derivative} of $F$.
  Goodwillie~\cites{go:ci,go:cii,go:ciii} has set up a ``calculus of
  functors'' in which the higher derivatives are spectra $D_nF$ with
  an action of $\Sg_n$.  (The slogan is that where the ordinary
  calculus of functions has a denominator of $n!$, the calculus of
  functors will take coinvariants under an action of $\Sg_n$.)  Even
  the derivatives of the identity functor are interesting; they fit in
  an intricate web of relationships with partition complexes,
  symmetric powers of the sphere spectrum, Steinberg modules, free Lie
  algebras and so on.  There are other versions of calculus for
  functors from other categories to spaces, with applications to
  embeddings of manifolds, for example.
 \item[(j)] A \emph{Moore spectrum} is a spectrum $X$ for which
  $\pi_n(X)=0$ when $n<0$ and $H_n(X)=0$ when $n>0$.  Let $\CH$ be the
  category of Moore spectra.  The functor $H_0\:\CH\to\Ab$ is then
  close to being an equivalence: for any $X,Y\in\CH$ there is a
  natural short exact sequence
  \[ \Ext(H_0(X),H_0(Y))/2 \xra{} [X,Y] \xra{}
      \Hom(\pi_0(X),\pi_0(Y)).
  \]
  Moreover, given any abelian group $A$ there is a Moore spectrum $SA$
  (unique up to non-canoncal isomorphism) with $H_0(SA)\simeq A$.
 \item[(k)] Let $C$ be an elliptic curve over a ring $k$.  (Number
  theorists are often interested in the case where $k$ is a small ring
  like $\Z$, but it is also useful to consider larger rings like
  $\Z[\tfrac{1}{6},c_4,c_6][(c_6^2-c_4^3)^{-1}]$ that have various
  universal properties in the theory of elliptic curves.)  From this
  we obtain a formal group $\widehat{C}$, which can be thought of as
  the part of $C$ infinitesimally close to zero.  It often happens
  that there is a spectrum $E$ that corresponds to $\widehat{C}$ under
  the standard dictionary relating formal groups to cohomology
  theories.  Spectra arising in this way are called \emph{elliptic
   spectra}~\cite{hoanst:esw}.  The details are usually adjusted so that
  $\pi_*(E)=k[u,u^{-1}]$ with $k=\pi_0(E)$ and $u\in\pi_2(E)$.  In
  many cases $E$ can be constructed using the Landweber Exact Functor
  Theorem, as in~(f).

  The spectrum $TMF$ (standing for \emph{topological modular forms})
  ``wants to be'' the universal example of an elliptic
  spectrum~\cite{dfhh:tmf}.  It is not in fact an elliptic spectrum,
  but it is close to being one, and it admits a canonical map to every
  elliptic spectrum.  If we let $MF_*$ denote the group of integral
  modular forms as defined by number theorists (graded so that forms
  of weight $k$ appear in degree $2k$) then we have
  \[ \pi_*(TMF)[\tfrac{1}{6}] =
       \Z[\tfrac{1}{6},c_4,c_6][(c_6^2-c_4^3)^{-1}] = 
         MF_*[\tfrac{1}{6}].
  \]
  The significance of the number $6$ here is that $2$ and $3$ are the
  only primes that can divide $|\Aut(C)|$, for any elliptic curve
  $C$.  If we do not invert $6$ then the homotopy groups $\pi_*(TMF)$
  are completely known, but different from $MF_*$ and too complex to
  describe here.

  There is a dense network of partially understood interactions
  between elliptic spectra, conformal nets~\cite{dohe:tmf}, vertex
  operator algebras, chiral differential operators~\cite{ch:wgv} and
  mathematical models of string theory~\cite{stte:sft}.  It seems
  likely that some central aspects of this picture remain to be
  discovered. 

  Under~(g) we discussed the algebraic $K$-theory spectrum $K(R)$
  associated to to a commutative ring $R$.  The construction can be
  generalised to define $K(R)$ when $R$ is a commutative ring
  spectrum; in particular, we can define $K(kU)$.  Rognes and
  Ausoni~\cite{auro:akt} have found evidence of a relationship between
  $K(kU)$ and $TMF$, but many features of this remain obscure.
 \item[(l)] Let $R$ be a strictly commutative ring spectrum, so for
  any space $X$ we have a ring $R^0(X)=[\Sgi X_+,R]$, and a group of
  units $R^0(X)^\tm$.  It can be shown that there is a spectrum
  $gl_1(R)$ such that $R^0(X)^\tm=[\Sgi X_+,gl_1(R)]$ for all $X$.  By
  taking $X=S^n$ we see that $\pi_n(gl_1(R))=\pi_n(R)$ for $n>0$, but
  $H_*(gl_1(R))$ is not closely related to $H_*(R)$, and many aspects
  of the topology of $gl_1(R)$ are mysterious, even in the case
  $R=S^0$.  For any $(-1)$-connected spectrum $T$ one can build a
  strictly commutative ring spectrum of homotopy type
  $\Sg^\infty(\Om^\infty T)_+$ (where $\Om^\infty T$ is the homotopy
  colimit of the spaces $\Om^nT(\R^n)$).  Ring maps from this to $R$
  then correspond to maps $T\to gl_1(R)$ of spectra.
 \item[(m)] The notion of Bousfield localisation~\cite{bo:lspe,ra:lrc}
  provides an important way to construct new spectra from old.  The
  simplest kinds of Bousfield localisations are the arithmetic ones:
  if $X$ is a spectrum and $p$ is a prime number then there are
  spectra and maps $X[\frac{1}{p}]\xla{i}X\xra{j}X_{(p)}$ such that
  $i$ induces an isomorphism 
  \[ \pi_*(X)[\tfrac{1}{p}] = \pi_*(X)\ot\Z[\tfrac{1}{p}]
      \to \pi_*(X[\tfrac{1}{p}])
  \]
  and $j$ induces an isomorphism 
  \[ \pi_*(X)_{(p)} = \pi_*(X)\ot\Z_{(p)} \to \pi_*(X_{(p)}). \]
  Similarly, there is a map $X\to X\Q$ inducing
  $\pi_*(X)\ot\Q\simeq\pi_*(X\Q)$.  These properties characterise
  $X[\frac{1}{p}]$, $X_{(p)}$ and $X\Q$ up to canonical homotopy
  equivalence.  Along similar lines, there is a $p$-adic completion
  map $X\to X^\wedge_p$.  If each homotopy group $\pi_k(X)$ is
  finitely generated, then we just have
  $\pi_k(X^\wedge_p)=\pi_k(X)^\wedge_p=\pi_k(X)\ot\Z_p$.  In the
  infinitely generated case the picture is more complicated, but still
  well-understood.  

  Next, for any spectrum $E$ there is a functor $L_E\:\CB\to\CB$ and a
  natural map $i\:X\to L_EX$ characterised as follows: the induced map
  $E_*i\:E_*X\to E_*L_EX$ is an isomorphism, and if $f\:X\to Y$ is
  such that $E_*f$ is an isomorphism, then there is a unique map
  $g\:Y\to L_EX$ with $gf=i$.  The spectrum $L_EX$ is called the
  \emph{Bousfield localisation} of $X$ with respect to $E$, and we say
  that $X$ is \emph{$E$-local} if the map $X\to L_EX$ is a homotopy
  equivalence.  The slogan is that the category $\CB_E$ of $E$-local
  spectra is the part of stable homotopy theory that is visible to
  $E$.  Apart from the arithmetic completions and localisations, the
  most important cases are the \emph{chromatic localisations}
  $L_{E(p,n)}$ and $L_{K(p,n)}$, which have been studied
  intensively~\cite{ra:lrc,host:mkl}.  It is a key point here that
  $E(p,n)$ and $K(p,n)$ have nontrivial homotopy groups in infinitely
  many negative degrees.  If the homotopy groups of $E$ and $X$ are
  bounded below, it turns out that $L_EX$ is just an arithmetic
  localisation of $X$, which is less interesting.
 \item[(n)] Under~(b) we mentioned the spectra $X(n)$ and
  $MU=X(\infty)$.  If we localise at a prime $p$ (as discussed in~(m))
  then these can be split into much smaller pieces.  For example, we
  have $\pi_*MU=\Z[x_1,x_2,\dotsc]$ with $|x_i|=2i$, but $MU_{(p)}$
  is equivalent to a coproduct of suspended copies of a spectrum
  called $BP$, where $\pi_*BP=\Z_{(p)}[v_1,v_2,\dotsc]$ with
  $|v_k|=2(p^k-1)$.  This is called the \emph{Brown-Peterson
   spectrum}.  There are also spectra called $T(n)$ that have the same
  relation to $X(p^n)$ as $BP$ does to $MU$.  This whole story is most
  naturally understood in terms of formal groups and $p$-typical
  curves.  The spectra $T(n)$ play an important role in a number of
  places, including the proof of the Nilpotence Theorem, the work of
  Mahowald, Ravenel and Shick on the Telescope
  Conjecture~\cite{marash:tls}, and Ravenel's methods for computation
  of stable homotopy groups of spheres~\cite{ra:ccs}.
 \item[(o)] The theory of surgery aims to understand compact smooth
  manifolds by cutting them into simpler pieces and reassembling the
  pieces~\cite{mami:css}.  A key ingredient is as follows: if we have
  an $n$-dimensional manifold with an embedded copy of
  $S^{i-1}\tm B^j$ (where $i+j=n+1$), then we can remove the interior
  to leave a manifold with boundary $S^{i-1}\tm S^{j-1}$, then glue on
  a copy of $B^i\tm S^{j-1}$ to obtain a new closed manifold $M'$.  If
  the original copy of $S^{i-1}\tm B^j$ is chosen appropriately, then
  the cohomology of $M'$ will be smaller than that of $M$.  By
  iterating this process, we hope to convert $M$ to $S^n$ (then we can
  reverse the steps to obtain a convenient description of $M$).  There
  are various obstructions to completing this process (and similar
  processes for related problems), and it turns out that these can be
  encoded as problems in stable homotopy theory.  Cobordism spectra
  such as $MSO$ play a role, as do the spectra $kO$ and $gl_1(S^0)$.
  When the dimension $n$ is even and $M$ is oriented, the
  multiplication map $H^{n/2}(M)\otimes H^{n/2}(M)\xra{}H^n(M)=\Z$
  gives a bilinear form on $H^{n/2}(M)$, which is symmetric or
  antisymmetric depending on the parity of $n/2$.  Because of this, it
  turns out that we need to consider a kind of $K$-theory of abelian
  groups equipped with a bilinear form.  This is known as $L$-theory.
  There are various different versions, depending on details that we
  have skipped over.  One of them gives a strictly commutative ring
  spectrum $L$ with $\pi_*(L)=\Z[x,y]/(2x,x^2)$, where $|x|=1$ and
  $|y|=4$ (so for $k\geq 0$ we have $\pi_{4k}(L)=\Z.y^k$ and
  $\pi_{4k+1}(L)=(\Z/2).y^kx$).  There is a canonical map
  $\sigma\:MSO\to L$ of strictly commutative ring
  spectra~\cite{lamc:mpq}.  For an oriented manifold $M$ of dimension
  $4k$, the cobordism class $[M]$ gives an element of $\pi_{4k}(MSO)$,
  so we must have $\sigma_*([M])=d.y^k$ for some integer $d$.  This
  integer is just the signature of the symmetric bilinear form on
  $H^{2k}(M;\R)$.  A similar description can be given in dimension
  $4k+1$. 
\end{itemize}

\end{document}